\documentclass[10pt,twoside]{article}
\usepackage{amsmath,amssymb,amsfonts}
\usepackage{latexsym}
\usepackage{txfonts}
\newtheorem{theorem}{Theorem}[section]

\newtheorem{lemma}{Lemma}[section]

\newtheorem{corollary}{Corollary}[section]

\makeatletter\def\theequation{\arabic{section}.\arabic{equation}}\makeatother
\begin{document}
\title{
\begin{flushleft}
\end{flushleft}
\vspace{0.5in}
{\bf\Large  A Quasilinear Elliptic Equation with Oscillating Nonlinearity}}
\author{{\bf\large Anderson L. A. de Araujo\vspace{1mm}}\\
{\it\small Departamento de Matem\'atica}\\ {\it\small Universidade Federal de Vi\c cosa},\\
{\it\small Av Peter Henry Rolfs, s/n, CEP 36570-900, Vi\c cosa, MG, Brasil}\\
{\it\small e-mail: anderson.araujo@ufv.br \vspace{1mm}}\vspace{1mm}\\
{\bf\large Rafael dos Reis Abreu}\\
{\it\small Departamento de Ci\^encias Exatas e Educa\c c\~ao do Centro de Blumenau}\\ {\it\small Universidade Federal de Santa Catarina},\\
{\it\small Rua Pomerode, 710, CEP 89065-300, Blumenau, SC, Brasil}\\
{\it\small e-mail: rafael.abreu@ufsc.br}
}\vspace{1mm}
\date{\small{Received XX April XX} \\
{\it\small Communicated ......}}
\maketitle
\begin{center}
{\bf\small Abstract}

\vspace{3mm}
\hspace{.05in}\parbox{4.5in}
{{\small We find a solution of the quasilinear elliptic equation
$-div\big(\phi(|\nabla u|^2)\nabla u\big)=f(x,u)$ in $\Omega$ with Dirichlet's boundary condition,
where $\Omega$ is a smooth bounded domain in $\mathbb{R}^N$, $N\geq 1$ and $f:\Omega\times \mathbb{R} \rightarrow \mathbb{R}$ is an unbounded continuous function with oscillatory behavior near the origin.}}
\end{center}
\section{\bf Introduction}
\def\theequation{1.\arabic{equation}}\makeatother
\setcounter{equation}{0}

We study the quasilinear elliptic second order equation
\begin{equation}\label{eq1}
\left\{
\begin{array}{rcll}
-div\big(\phi(|\nabla u|^2)\nabla u\big)  & =  &f(x,u)\,\, \mbox{ in } \,\,  \Omega, \\
u & =  &0 \,\, \mbox{ on }\,\, \partial \Omega,
\end{array}
\right.
\end{equation}
where $\Omega$ is a smooth bounded domain in $\mathbb{R}^N$, $N\geq 1$, $- div\big(\phi(|\nabla u|^2) \nabla u \big)$ is the $\phi$- Laplacian operator and $f:\Omega \times [0, + \infty) \rightarrow \mathbb{R}$ is an unbounded continuous function with oscillatory behavior near the origin in the sense that follows. There are two sequences of real numbers $\beta_n >0$ and $\gamma_n > 0$ such that
\begin{equation}\label{desigualdade}
\beta_n \leq \gamma_n \ , \ \forall n \in \mathbb{N} ;
\end{equation}
\begin{equation}\label{eq2.0}
\beta_n \rightarrow 0^+ \ \ \textrm{and} \ \ \gamma_n \rightarrow 0^+;
\end{equation}
\begin{equation}\label{eq2}
f(x,\gamma_n)\leq 0 \,\,\, \mbox{in}\,\,\, \Omega \ \ \ \forall n \in \mathbb{N};	
\end{equation}
and
\begin{equation}\label{eq3}
\liminf_n\frac{f(x,\beta_n\varphi_1)}{\beta_n}=+\infty \,\,\,\mbox{uniformly in}\,\,\,\Omega,	
\end{equation}
where $\varphi_1$ is an eigenfunction of the Laplacian operator in $\Omega$ corresponding to the first eigenvalue $\lambda_1 > 0$, that is, $\varphi_1 \in H_0^1(\Omega) \cap C^{\infty}(\overline{\Omega})$ satisfies $\varphi_1 > 0$ in $\Omega$ and
\begin{eqnarray*}\label{eq5}
\left\{ \begin{array}{ccc} -\Delta \varphi_1=\lambda_1\varphi_1 \,\,\mbox{in}\,\,\Omega, \\ \varphi_1 = 0 \ \ \textrm{on} \ \ \partial \Omega. \end{array} \right.
\end{eqnarray*}
By normalization, we can suppose $0< \varphi_1 \leq 1$ in $\Omega$.

An example of a function which satisfies the above hypotheses is
\begin{eqnarray}\label{f}
f(x,t)=\left\{
\begin{array}{lcl}
\displaystyle  -\lambda_1 \sqrt{\frac{t}{\varphi_1(x)}} \sin\left(\frac{\sqrt{\varphi_1(x)}}{\sqrt{t}}\right) &,& t>0, \\
0&,& t=0.
\end{array}
\right.
\end{eqnarray}
For each $n \in \mathbb{N}$, let us define 
\begin{equation}\label{betan}
	\beta_n= \left( \frac{1}{2n\pi + \frac{3\pi}{2}} \right)^2 
\end{equation}
and
\begin{equation*}\label{gamman}
 \gamma_n = \left( \frac{1}{2n\pi + \frac{\pi}{4}} \right)^2.
\end{equation*}
Hence
\begin{eqnarray*}
\beta_n \rightarrow 0^+ \ \ \textrm{and} \ \ \gamma_n \rightarrow 0^+;
\end{eqnarray*}
\begin{eqnarray*}
\beta_n \leq \gamma_n \ \ , \ \ \forall n \in \mathbb{N};
\end{eqnarray*}
\begin{eqnarray*}
f(x,\gamma_n)= -\lambda_1 \frac{1}{\sqrt{\varphi_1(x)} \left( 2n \pi + \frac{\pi}{4}\right)} \sin\left(\sqrt{\varphi_1(x)} \left( 2n\pi + \frac{\pi}{4}\right)\right)\leq 0 \,\,\, \mbox{in}\,\,\, \Omega \ \ \ \forall n \in \mathbb{N}
\end{eqnarray*}
and
\[\frac{f(x,\beta_n\varphi_1)}{\beta_n}= \lambda_1 \left(2n\pi + \frac{3\pi}{2}\right).\]
Therefore, $f$, $\beta_n$ and $\gamma_n$ just defined satisfy (\ref{desigualdade}), (\ref{eq2.0}), (\ref{eq2}) and (\ref{eq3}).

We suppose that the function $\phi:[0,+\infty[ \to [0,+\infty[$ is of class $C^1$ and satisfies the following ellipticity and growth conditions of Leray-Lions type (see \cite{Tol}): there are constants $\gamma, \Gamma >0$, $\kappa \in [0,1]$ and $p \in [2,+\infty[$ such that, for every $s>0$,
\begin{equation}\label{ellipt}
\gamma\,(\kappa + s)^{p-2}\leq \phi(s^2)\leq \Gamma\,(\kappa + s)^{p-2}	
\end{equation}
and
\begin{equation}\label{LL}
(\gamma-\frac{1}{2})\phi(s)\leq \phi'(s)s\leq \Gamma\,\phi(s).	
\end{equation}
Hence, we can define, for $s\geq 0$, the function $\Phi(s)=\int_0^s \phi(\xi)d\xi$.

We also suppose that there is $t_0>0$ such that, for every $w \in H_0^1(\Omega)$ satisfying $w \geq 0$ a.e. in $\Omega$ and for $0<\beta_n<t_0$, we have
\begin{equation}\label{In1}
\int_{\Omega}\phi(|\nabla\,\beta_n\varphi_1|^2)\nabla\,\beta_n\varphi_1\nabla\,wdx \leq \int_{\Omega}f(x,\beta_n\varphi_1)wdx.
\end{equation}
Note that, with the assumptions on $\phi$, the function $f$ and the sequence $(\beta_n)$ defined in (\ref{f}) and (\ref{betan}), respectively, satisfy \eqref{In1}. In fact, by integration by parts, we have
\[
\begin{array}{rcl}
\displaystyle \int_{\Omega}\phi(|\nabla\,\beta_n\varphi_1|^2)\nabla\,\beta_n\varphi_1\nabla\,wdx &=&\displaystyle -\int_{\Omega}div\left(\phi(|\nabla\,\beta_n\varphi_1|^2)\nabla\,\beta_n\varphi_1\right)wdx\\
&=&\displaystyle -2\beta_n^3\int_{\Omega}\left(\phi'(|\nabla\,\beta_n\varphi_1|^2)\sum_{i=1}^N\partial_{x_i}\varphi_1\left\langle \nabla\varphi_1,\nabla(\partial_{x_i}\varphi_1)\right\rangle\right)wdx\\
&+&\displaystyle \lambda_1\beta_n\int_{\Omega}\phi(|\nabla\,\beta_n\varphi_1|^2)\varphi_1wdx.
\end{array}
\]
Then, in this case, a sufficient condition to occur \eqref{In1} is 
\[
\begin{array}{c}
\displaystyle -2\beta_n^3\int_{\Omega}\left(\phi'(|\nabla\,\beta_n\varphi_1|^2)\sum_{i=1}^N\partial_{x_i}\varphi_1\left\langle \nabla\varphi_1,\nabla(\partial_{x_i}\varphi_1)\right\rangle\right)wdx + \displaystyle \lambda_1\beta_n\int_{\Omega}\phi(|\nabla\,\beta_n\varphi_1|^2)\varphi_1wdx\\ 
\displaystyle \leq \int_{\Omega}f(x,\beta_n\varphi_1)wdx,
\end{array}
\]
which is equivalent to
\[
\int_{\Omega}\left[ \lambda_1 \left(2n\pi + \frac{3\pi}{2}\right) - \lambda_1\,\phi(|\nabla\,\beta_n\varphi_1|^2)\varphi_1 + 2\beta_n^2\left(\phi'(|\nabla\,\beta_n\varphi_1|^2)\sum_{i=1}^N\partial_{x_i}\varphi_1\left\langle \nabla\varphi_1,\nabla(\partial_{x_i}\varphi_1)\right\rangle\right)\right] wdx\geq 0
\]
since
\begin{eqnarray*}
f(x,\beta_n \varphi_1) = \lambda_1 \beta_n \left( 2n\pi + \frac{3 \pi}{2} \right).
\end{eqnarray*}
Provided that $\varphi_1 \in C^{\infty}(\overline{\Omega})$, it follows from \eqref{eq3}, \eqref{ellipt} and \eqref{LL} the existence of $n_0>0$ such that, for $n\geq n_0$, we have
\[
\lambda_1 \left(2n\pi + \frac{3\pi}{2}\right) - \lambda_1\,\phi(|\nabla\,\beta_n\varphi_1|^2)\varphi_1 + 2\beta_n^2\left(\phi'(|\nabla\,\beta_n\varphi_1|^2)\sum_{i=1}^N\partial_{x_i}\varphi_1\left\langle \nabla\varphi_1,\nabla(\partial_{x_i}\varphi_1)\right\rangle\right)\geq 0
\]
and, since $w\geq 0$, we verify the assumptions \eqref{In1}.

Quasilinear problems involving the $\phi$-Laplacian operator have been studied by some authors. In \cite{Njoku}, the authors studied the following quasilinear elliptic problem in a ball involving an oscillatory nonlinearity:
\begin{eqnarray*}
\left\{
\begin{array}{ccc}
-div\big(\phi(|\nabla u|^2)\nabla u\big)  & =  &f(|x|,u)\,\, \mbox{ in } \,\,  B_R, \\
u & =  &0 \,\, \mbox{ on }\,\, \partial B_R,
\end{array}
\right.
\end{eqnarray*}
where $B_R$ denotes the open ball at the origin with radius $R$ in $\mathbb{R}^N$, $\phi:[0, +\infty[ \rightarrow [0, + \infty[$ is continuous and satisfies
\begin{eqnarray*}
\liminf_{s \rightarrow + \infty} \frac{\phi((\sigma s)^2)}{\phi(s^2)} > 1 \ \ \ \textrm{and} \ \ \ \limsup_{s \rightarrow + \infty} \frac{\phi((\sigma s)^2)}{\phi(s^2)} < + \infty
\end{eqnarray*}
for all $\sigma > 1$ and $f: [0, R] \times \mathbb{R} \rightarrow \mathbb{R}$ is a continuous function such that $s \mapsto \int^{s}_{0} f(|x|, \xi) d \xi$ has an oscillatory behaviour at $+\infty$. About quasilinear problems involving the $\phi$-Laplacian operator, we can also refer \cite{Carvalho, Schmitt, OZ, Wang}.

It is also worth noting that problems with oscillatory nonlinearities have been of interest to many authors; see, for example, \cite{AM, Kristaly, Malin, Mawhin, Ob, OZ}.

In this work, we are concerned with the existence of a weak positive solution of (\ref{eq1}). By a weak positive solution, we mean a function $u \in W^{1,p}_0(\Omega)\cap L^{\infty}(\Omega)$ satisfying $u > 0$ a.e. in $\Omega$ and
\[
\int_{\Omega}\phi(|\nabla u|^2)\nabla u\nabla wdx	=\int_{\Omega}f(x,u)wdx,
\]
for every $w \in W_0^{1,p}(\Omega)$, where $p$ is the exponent which appears in \eqref{ellipt}.	We emphasize that we complement the result obtained in \cite{Njoku} in the sense that we prove existence of a solution when $\Omega$ is not necessarily a ball and $f$ has an oscillatory behavior near the origin. Furthermore, our hypothesis on the $\phi$-Laplacian operator are different of those used there.

We state the main result.
\begin{theorem}\label{T1}
Let $\Omega$ be a smooth bounded domain in $\mathbb{R}^N$, $N\geq 1$. Assume that $\phi: [0,+\infty[\rightarrow [0,+\infty[$ satisfies (\ref{ellipt}), (\ref{LL}), $f:\Omega \times \mathbb{R} \rightarrow \mathbb{R}$ is a continuous function satisfying (\ref{desigualdade}), (\ref{eq2.0}), (\ref{eq2}) and (\ref{eq3}) and (\ref{In1}) is satisfied. Then, there is $n_0 \in \mathbb{N}$ such that, for every $n\geq n_0$, the problem (\ref{eq1}) has at least one positive solution $u_n$. Furthermore, this solution satisfies $\|u_n\|_{C^1(\overline{\Omega})} \to 0$ as $n\to +\infty$.
\end{theorem}

The class of quasilinear differential operators which we can consider includes in particular the $p$-Laplacian operator $\Delta_p$, with $2\leq p<+\infty$; in this case, $\phi(s)=s^{\frac{p}{2}-1}$ and $\Phi(s^2)=\frac{2}{p}|s|^p$. The Theorem \ref{T1} yields the following corollary.
\begin{corollary}\label{C1}
Assume the same assumptions in Theorem \ref{T1}. Then, the same conclusions of Theorem \ref{T1} holds for the problem
\[
\left\{ \begin{array}{rcll}
-\Delta_pu  & =  &f(x,u)\,\, \mbox{ in } \,\,  \Omega, \\
u & =  &0 \,\, \mbox{ on }\,\, \partial \Omega.
\end{array} \right.
\]

\end{corollary}

\section{\bf Proof of the main result}

We say that a function $\beta \in C^{1}(\overline{\Omega})$ is a \textsl{supersolution} of (\ref{eq1}) if
\[
\beta(x)\geq 0\,\,\mbox{on}\,\,\partial\Omega
\]
and
\[
\int_{\Omega}\phi(|\nabla \beta|^2)\nabla \beta \nabla wdx \geq \int_{\Omega}f(x,\beta)wdx
\]
for every $w \in W_0^{1,p}(\Omega)$, with $w \geq 0$ a.e. in $\Omega$.

Similarly, $\alpha \in C^{1}(\overline{\Omega})$ is a \textsl{subsolution} of (\ref{eq1}) if
\[
\alpha(x)\leq 0\,\,\mbox{on}\,\,\partial\Omega
\]
and
\[
\int_{\Omega}\phi(|\nabla \alpha |^2)\nabla \alpha \nabla wdx \leq \int_{\Omega} f(x,\alpha)wdx
\]
for every $w \in W_0^{1,p}(\Omega)$, with $w \geq 0$ a.e. in $\Omega$.


We also introduce the functional $J: W_0^{1,p}(\Omega)\cap L^{\infty}(\Omega) \rightarrow \mathbb{R}$ defined by
\begin{equation*}\label{eq8}
J(u)=\frac{1}{2}\int_{\Omega}\Phi(|\nabla u|^2)dx - \int_{\Omega}F(x,u)dx,
\end{equation*}
where $\Phi(v)=\int_0^v \phi(s)ds$ and $F(x,u)=\int_0^u f(x,s)ds$. 

The following existence and regularity result can be easily derived by combining \cite[Lemma 2.1]{OZ} and \cite{Lie}, see also \cite[Lemma 3.1]{HO}.

\begin{lemma}\label{lem1}
Let $\Omega$ be as in Theorem \ref{T1}, $\phi: [0,+\infty[\rightarrow [0,+\infty[$ satisfying (\ref{ellipt}), (\ref{LL}) and $f:\Omega\times[0, + \infty[\rightarrow~\mathbb{R}$ is a continuous function. Assume that there exist a subsolution
$\alpha$ and a supersolution $\beta$ of (\ref{eq1}) with $\alpha(x)\leq \beta(x)$ in $\Omega$. Then problem (\ref{eq1}) has at least one solution $u$ such that $\alpha(x)\leq u(x)\leq \beta(x)$ in $\Omega$ and
\[J(u)=\min\{J(v)| v \in W_0^{1,p}(\Omega)\,\,\mbox{and}\,\,\alpha(x)\leq v(x)\leq \beta(x)\,\,\mbox{a.e.}\,\,\Omega\}.\]
Moreover, there are constants $\sigma>0$ and $C>0$ continuously depending only on $\|\alpha\|_{\infty}$, $\|\beta\|_{\infty}$ and $\max\{|f(x,u)|| x \in \overline{\Omega}, \alpha(x)\leq u\leq \beta(x)\}$ such that, for every solution $u$ of (\ref{eq1}) with  $\alpha(x)\leq u(x)\leq \beta(x)$ in $\Omega$, we have $u \in C^{1,\sigma}(\overline{\Omega})$ and
\[\|u\|_{C^{1,\sigma}(\overline{\Omega})}\leq C.\]
\end{lemma}

\subsection{Proof of Theorem \ref{T1}}

Notice that for all $n \in \mathbb{N}$, $\gamma_n$ is a supersolution of (\ref{eq1}). Indeed, $\gamma_n \in C^1(\overline{\Omega})$, $\gamma_n \geq 0$ in $\partial\Omega$ and, by (\ref{eq2}),
\[
\int_{\Omega}\phi(|\nabla\gamma_n|^2)\nabla\gamma_n\nabla\,wdx =0\geq \int_{\Omega}f(x,\gamma_n)wdx
\]
for every $w \in H_0^1(\Omega)$, with $w(x) \geq 0$ a.e. in $\Omega$. Therefore, $\gamma_n$ is a supersolution of (\ref{eq1}).

On the other hand, notice that, for $n \in \mathbb{N}$ large enough, $0<\beta_n<t_0$ and it follows from (\ref{In1}) that $\beta_n\varphi_1$ is a subsolution of (\ref{eq1}).

From (\ref{desigualdade}) and the hypothesis of normalization of $\varphi_1$, it follows that
\begin{eqnarray*}
\beta_n \varphi_1(x)\leq \beta_n \leq \gamma_n \ \ , \ \ \forall x \in \Omega
\end{eqnarray*}
for all $n \in \mathbb{N}$ large enough. Then, by Lemma \ref{lem1}, we have that there exists $n_0 \in \mathbb{N}$ such that the problem (\ref{eq1}) has, for each $n \geq n_0$, at least one solution $u_n$ such that
\begin{eqnarray}\label{eq9}
	\beta_n\varphi_1\leq u_n\leq \gamma_n\,\,\,\mbox{in}\,\,\Omega,
\end{eqnarray}
\[J(u_n)=\min\{J(v)| v \in W_0^{1,p}(\Omega)\,\,\mbox{and}\,\,\beta_n\varphi_1(x)\leq v(x)\leq \gamma_n\,\,\mbox{a.e.}\,\,\Omega\}\]
and
\begin{eqnarray}\label{eq10}
\|u_n\|_{C^{1,\tau}(\overline{\Omega})}\leq C,
\end{eqnarray}
where the constants $1\geq \tau>0$ and $C>0$ is independent of $n$.
%
By Arzel\`a-Ascoli Theorem and (\ref{eq10}), there is a subsequence $(u_{n_k})$ such that
\[
u_{n_k} \rightarrow u\,\,\,\mbox{in}\,\,\, C^1(\overline{\Omega}),
\]
with $u=0$ em $\partial\Omega$. But letting $k \rightarrow + \infty$, we obtain $\beta_{n_k}\rightarrow 0$ and $\gamma_{n_k}\rightarrow 0$; then, by (\ref{eq9}) , we conclude that $u_{n_k} \rightarrow 0$ uniformly in $\Omega$ and therefore $u=0$.

Hence, by (\ref{eq10}), the sequence $(u_n)$ has the property that every subsequence has a subsequence which converges to zero in $C^1(\overline{\Omega})$. In conclusion,
\[u_{n} \rightarrow 0\,\,\,\mbox{in}\,\,\, C^1(\overline{\Omega}).\]
From the above results, there exists $n_0 \in \mathbb{N}$ such that, for $n\geq n_0$, we have
\[	\beta_n\varphi_1\leq u_n\leq \gamma_n\,\,\,\mbox{in}\,\,\Omega\]
and
\[|\nabla u_n(x)|<1\,\,\,\mbox{in}\,\,\Omega.\]
Consequently, for $n\geq n_0$, $u_n$ is a solution of (\ref{eq1}) and
\[\lim_{n\rightarrow \infty}\|u_n\|_{C^{1}(\overline{\Omega})}=0.\]

\end{document}